\documentclass{amsart}
\usepackage{amssymb}
\usepackage{amsfonts}

\newtheorem{theorem}{Theorem}
\theoremstyle{plain}

\newtheorem{corollary}{Corollary}

\newtheorem{definition}{Definition}

\numberwithin{equation}{section}
\input{tcilatex}

\begin{document}
\title[Fractional type higher order commutators]{The boundedness of a class
of fractional type rough higher order commutators {on vanishing generalized
weighted Morrey spaces}}
\author{FER\.{I}T G\"{U}RB\"{U}Z}
\address{HAKKARI UNIVERSITY, FACULTY OF EDUCATION, DEPARTMENT OF MATHEMATICS
EDUCATION, HAKKARI 30000, TURKEY }
\email{feritgurbuz@hakkari.edu.tr}
\urladdr{}
\thanks{}
\curraddr{ }
\urladdr{}
\thanks{}
\date{}
\dedicatory{}
\thanks{}

\begin{abstract}
This paper includes new bounds concepting the {vanishing generalized
weighted Morrey space. In this sense,} it is outlined improved bounds about
the a class of fractional type rough higher order commutators {on vanishing
generalized weighted Morrey spaces.}
\end{abstract}

\maketitle

\section{Introduction}

Let $\Omega \in L_{s}(S^{n-1})$, $1<s\leq \infty $. $\Omega $ is the
function defined on ${\mathbb{R}^{n}}\setminus \{0\}$ satisfying the
homogeneous of degree zero condition, that is,%
\begin{equation}
\Omega (\lambda x)=\Omega (x)\text{~for any~}\lambda >0\text{, }x\in {%
\mathbb{R}^{n}}\setminus \{0\}  \label{f}
\end{equation}%
and the integral zero property (=the vanishing moment condition) over the
unit sphere $S^{n-1}$, that is, 
\begin{equation}
\int\limits_{S^{n-1}}\Omega (x^{\prime })d\sigma (x^{\prime })=0,  \label{0*}
\end{equation}%
where $x^{\prime }=\frac{x}{|x|}$ for any $x\neq 0$.

In this paper we consider the following higher order ($=k$-th order)
commutator operators of rough fractional integral and maximal operators,%
\begin{eqnarray}
T_{\Omega ,\alpha }^{A,k}f(x) &=&T_{\Omega ,\alpha }\left( \left( A\left(
x\right) -A\left( \cdot \right) \right) ^{k}f\left( \cdot \right) \right)
\left( x\right) ,\qquad k=0,1,2,\ldots ,  \notag \\
&=&\int\limits_{{\mathbb{R}^{n}}}\frac{\Omega (x-y)}{|x-y|^{n-\alpha }}%
\left( A\left( x\right) -A\left( y\right) \right) ^{k}f(y)dy  \label{f1}
\end{eqnarray}%
and 
\begin{eqnarray}
M_{\Omega ,\alpha }^{A,k}f(x) &=&M_{\Omega ,\alpha }\left( \left( A\left(
x\right) -A\left( \cdot \right) \right) ^{k}f\left( \cdot \right) \right)
\left( x\right) ,\qquad k=0,1,2,\ldots ,  \notag \\
&=&\sup_{r>0}\frac{1}{r^{n-\alpha }}\int\limits_{|x-y|<r}\left\vert \Omega
\left( x-y\right) \right\vert \left\vert A\left( x\right) -A\left( y\right)
\right\vert ^{k}\left\vert f(y)\right\vert dy,  \label{f2}
\end{eqnarray}%
as long as the integrals above make sense, where rough fractional integral
operator $T_{\Omega ,\alpha }$ and rough fractional maximal operator $%
M_{\Omega ,\alpha }$ are defined by%
\begin{equation*}
T_{\Omega ,\alpha }f(x)=\int\limits_{{\mathbb{R}^{n}}}\frac{\Omega (x-y)}{%
|x-y|^{n-\alpha }}f(y)dy\qquad 0<\alpha <n
\end{equation*}%
and 
\begin{equation*}
M_{\Omega ,\alpha }f(x)=\sup_{r>0}\frac{1}{r^{n-\alpha }}\dint%
\limits_{|x-y|<r}\left\vert \Omega (x-y)\right\vert \left\vert
f(y)\right\vert dy\qquad 0<\alpha <n.
\end{equation*}

For $k=1$ above, $T_{\Omega ,\alpha }^{A,k}$ and $M_{\Omega ,\alpha }^{A,k}$
are obviously reduced to the rough commutator operators of $T_{\Omega
,\alpha }$ and $M_{\Omega ,\alpha }$, respectively:%
\begin{eqnarray*}
\left[ A,T_{\Omega ,\alpha }\right] f\left( x\right) &=&A\left( x\right)
T_{\Omega ,\alpha }f\left( x\right) -T_{\Omega ,\alpha }\left( Af\right)
\left( x\right) \\
&=&\int\limits_{{\mathbb{R}^{n}}}\frac{\Omega (x-y)}{|x-y|^{n-\alpha }}%
\left( A\left( x\right) -A\left( y\right) \right) f(y)dy
\end{eqnarray*}%
and 
\begin{eqnarray*}
\left[ A,M_{\Omega ,\alpha }\right] f\left( x\right) &=&A\left( x\right)
M_{\Omega ,\alpha }f\left( x\right) -M_{\Omega ,\alpha }\left( Af\right)
\left( x\right) \\
&=&\sup_{r>0}\frac{1}{r^{n-\alpha }}\int\limits_{|x-y|<r}\left\vert \Omega
\left( x-y\right) \right\vert \left\vert A\left( x\right) -A\left( y\right)
\right\vert \left\vert f(y)\right\vert dy.
\end{eqnarray*}%
Moreover, $T_{\Omega ,\alpha }^{A,k}$ and $M_{\Omega ,\alpha }^{A,k}$ are
trivial generalizations of the above commutators, respectively.

Here and henceforth, $F\approx G$ means $F\gtrsim G\gtrsim F$; while $%
F\gtrsim G$ means $F\geq CG$ for a constant $C>0$; and also $C$ stands for a
positive constant that can change its value in each statement without
explicit mention.

Now, let us list some definitions that we need in the proof of following
Theorem \ref{teo3}:

\begin{definition}
$\left( \text{\textbf{Bounded Mean Oscillation (BMO)}}\right) $ We denote
the mean value of $f$ on $B=B(x,r)\subset {\mathbb{R}^{n}}$ by%
\begin{equation*}
f_{B}=M\left( f,B\right) =M\left( f,x,r\right) =\frac{1}{|B|}%
\dint\limits_{B}f(y)dy,
\end{equation*}%
and the mean oscillation of $f$ on $B=B(x,r)$ by%
\begin{equation*}
MO\left( f,B\right) =MO\left( f,x,r\right) =\frac{1}{|B|}\dint%
\limits_{B}|f(y)-f_{B}|dy.
\end{equation*}%
We also define for a non-negative function $\phi $ on ${\mathbb{R}^{n}}$%
\begin{equation*}
MO_{\phi }\left( f,B\right) =MO_{\phi }\left( f,x,r\right) =\frac{1}{\phi
\left( |B|\right) |B|}\dint\limits_{B}|f(y)-f_{B}|dy.
\end{equation*}%
Now we define%
\begin{equation*}
BMO_{\phi }=\left\{ f\in L_{1}^{loc}({\mathbb{R}^{n}}):\sup_{B}MO_{\phi
}\left( f,B\right) <\infty \right\}
\end{equation*}%
and%
\begin{equation*}
\left\Vert f\right\Vert _{BMO_{\phi }}=\sup_{B}MO_{\phi }\left( f,B\right) .
\end{equation*}%
The most important of these spaces occurs when $\phi =1$, in which case $%
BMO_{\phi }=BMO$.
\end{definition}

\begin{definition}
\cite{Gurbuz2, Gurbuz4}$\left( \text{\textbf{Weighted Lebesgue space}}%
\right) $ Let $1\leq p\leq \infty $ and given a weight $w\left( x\right) \in
A_{p}\left( {{\mathbb{R}^{n}}}\right) $, we shall define weighted Lebesgue
spaces as 
\begin{eqnarray*}
L_{p}(w) &\equiv &L_{p}({{\mathbb{R}^{n}}},w)=\left\{ f:\Vert f\Vert
_{L_{p,w}}=\left( \dint\limits_{{{\mathbb{R}^{n}}}}|f(x)|^{p}w(x)dx\right) ^{%
\frac{1}{p}}<\infty \right\} ,\qquad 1\leq p<\infty . \\
L_{\infty ,w} &\equiv &L_{\infty }({{\mathbb{R}^{n}}},w)=\left\{ f:\Vert
f\Vert _{L_{\infty ,w}}=\limfunc{esssup}\limits_{x\in {\mathbb{R}^{n}}%
}|f(x)|w(x)<\infty \right\} .
\end{eqnarray*}
\end{definition}

\bigskip Here and later, we refer to $A_{p}$ as the the Muckenhoupt classes.
That is, $w\left( x\right) \in A_{p}\left( {{\mathbb{R}^{n}}}\right) $ for
some $1<p<\infty $ if $\left( \frac{1}{|B|}\dint\limits_{B}w(y)dy\right)
\left( \frac{1}{|B|}\dint\limits_{B}w(y)^{-\frac{1}{p-1}}dy\right)
^{p-1}\leq C$ for all balls $B$ (see \cite{Gurbuz2} for more details).

Now, let us consider the Muckenhoupt-Wheeden class $A\left( p,q\right) $ in 
\cite{Muckenhoupt}. One says that $w\left( x\right) \in A\left( p,q\right) $
for $1<p<q<\infty $ if and only if 
\begin{equation}
\lbrack w]_{A\left( p,q\right) }:=\sup\limits_{B}\left(
|B|^{-1}\dint\limits_{B}w(x)^{q}dx\right) ^{\frac{1}{q}}\left(
|B|^{-1}\dint\limits_{B}w(x)^{-p^{\prime }}dx\right) ^{\frac{1}{p^{\prime }}%
}<\infty ,  \label{13}
\end{equation}%
where the supremum is taken over all the balls $B$. Note that, by H\"{o}%
lder's inequality, for all balls $B$ we have 
\begin{equation*}
\lbrack w]_{A\left( p,q\right) }\geq \lbrack w]_{A\left( p,q\right)
(B)}=|B|^{\frac{1}{p}-\frac{1}{q}-1}\Vert w\Vert _{L_{q}(B)}\Vert
w^{-1}\Vert _{L_{p^{\prime }}(B)}\geq 1.
\end{equation*}

By (\ref{13}), we have%
\begin{equation*}
\left( \dint\limits_{B}w(x)^{q}dx\right) ^{\frac{1}{q}}\left(
\dint\limits_{B}w(x)^{-p^{\prime }}dx\right) ^{\frac{1}{p^{\prime }}%
}\lesssim \left\vert B\right\vert ^{\frac{1}{q}+\frac{1}{p^{\prime }}}.
\end{equation*}%
On the other hand, let $\mu \left( x\right) =w\left( x\right) ^{s^{\prime }}$%
, $\tilde{p}=\frac{p}{s^{\prime }}$ and $\tilde{q}=\frac{q}{s^{\prime }}$.
If $w\left( x\right) ^{s^{\prime }}\in A\left( \frac{p}{s^{\prime }},\frac{q%
}{s^{\prime }}\right) $, then we get $\mu \left( x\right) \in A\left( \tilde{%
p},\tilde{q}\right) $.

Now, we introduce some spaces which play important roles in PDE. Except the
weighted Lebesgue space $L_{p}(w)$, the weighted Morrey space $L_{p,\kappa
}(w)$, which is a natural generalization of $L_{p}(w)$ is another important
function space. Then, the definition of generalized weighted Morrey spaces $%
M_{p,\varphi }\left( w\right) $ which could be viewed as extension of $%
L_{p,\kappa }(w)$ has been given as follows:

For $1\leq p<\infty $, positive measurable function $\varphi (x,r)$ on ${%
\mathbb{R}^{n}}\times (0,\infty )$ and nonnegative measurable function $w$
on ${\mathbb{R}^{n}}$, $f\in M_{p,\varphi }(w)\equiv M_{p,\varphi }({\mathbb{%
R}^{n}},w)$ if $f\in L_{p,w}^{loc}({\mathbb{R}^{n}})$ and 
\begin{equation*}
\Vert f\Vert _{M_{p,\varphi }(w)}=\sup\limits_{x\in {\mathbb{R}^{n}},r>0}%
\frac{1}{\varphi (x,r)}\Vert f\Vert _{L_{p}(B(x,r),w)}<\infty .
\end{equation*}%
is finite. Note that for $\varphi (x,r)\equiv w(B(x,r))^{\frac{\kappa }{p}}$%
, $0<\kappa <1$ and $\varphi (x,r)\equiv 1$, we have $M_{p,\varphi
}(w)=L_{p,\kappa }(w)$ and $M_{p,\varphi }(w)=L_{p}(w)$, respectively.
Moreover, G\"{u}rb\"{u}z \cite{Gurbuz3} proved that the operators $T_{\Omega
,\alpha }^{A,k}$ and $M_{\Omega ,\alpha }^{A,k}$ are bounded from one
generalized weighted Morrey space $M_{p,\varphi _{1}}\left( w^{p},{\mathbb{R}%
^{n}}\right) $ to another $M_{q,\varphi _{2}}\left( w^{q},{\mathbb{R}^{n}}%
\right) $.

The following definition was introduced by G\"{u}rb\"{u}z \cite{Gurbuz5}.

\begin{definition}
\label{definition1}\textbf{(Vanishing generalized weighted Morrey spaces) }%
For\textbf{\ }$1\leq p<\infty $, $\varphi (x,r)$ is a positive measurable
function on ${\mathbb{R}^{n}}\times (0,\infty )$ and nonnegative measurable
function $w$ on ${\mathbb{R}^{n}}$, $f\in VM_{p,\varphi }\left( w\right)
\equiv VM_{p,\varphi }({\mathbb{R}^{n},w})$ if $f\in L_{p,w}^{loc}({\mathbb{R%
}^{n}})$ and%
\begin{equation}
\lim\limits_{r\rightarrow 0}\sup\limits_{x\in {\mathbb{R}^{n}}}\frac{1}{%
\varphi (x,r)}\Vert f\Vert _{L_{p}(B(x,r),w)}=0.  \label{1*}
\end{equation}
\end{definition}

Inherently, it is appropriate to impose on $\varphi (x,t)$ with the
following circumstances:

\begin{equation}
\lim_{t\rightarrow 0}\sup\limits_{x\in {\mathbb{R}^{n}}}\frac{\left(
w(B(x,t))\right) ^{^{\frac{1}{p}}}}{\varphi (x,t)}=0,  \label{2}
\end{equation}%
and%
\begin{equation}
\inf_{t>1}\sup\limits_{x\in {\mathbb{R}^{n}}}\frac{\left( w(B(x,t))\right)
^{^{\frac{1}{p}}}}{\varphi (x,t)}>0.  \label{3}
\end{equation}

From (\ref{2}) and (\ref{3}), we easily know that the bounded functions with
compact support belong to $VM_{p,\varphi }\left( w\right) $. On the other
hand, the space $VM_{p,\varphi }(w)$ is Banach space with respect to the
following finite quasi-norm%
\begin{equation*}
\Vert f\Vert _{VM_{p,\varphi }(w)}=\sup\limits_{x\in {\mathbb{R}^{n}},r>0}%
\frac{1}{\varphi (x,r)}\Vert f\Vert _{L_{p}(B(x,r),w)},
\end{equation*}%
such that%
\begin{equation*}
\lim\limits_{r\rightarrow 0}\sup\limits_{x\in {\mathbb{R}^{n}}}\frac{1}{%
\varphi (x,r)}\Vert f\Vert _{L_{p}(B(x,r),w)}=0,
\end{equation*}%
we omit the details. Moreover, we have the following embeddings:%
\begin{equation*}
VM_{p,\varphi }\left( w\right) \subset M_{p,\varphi }\left( w\right) ,\qquad
\Vert f\Vert _{M_{p,\varphi }\left( w\right) }\leq \Vert f\Vert
_{VM_{p,\varphi }\left( w\right) }.
\end{equation*}%
Henceforth, we denote by $\varphi \in \mathcal{B}\left( w\right) $ if $%
\varphi (x,r)$ is a positive measurable function on ${\mathbb{R}^{n}}\times
(0,\infty )$ and positive for all $(x,r)\in {\mathbb{R}^{n}}\times (0,\infty
)$ and satisfies (\ref{2}) and (\ref{3}).

Inspired of \cite{Gurbuz3}, the aim of the present paper is to study the
boundedness of the operators $T_{\Omega ,\alpha }^{A,k}$ and $M_{\Omega
,\alpha }^{A,k}${\ generated by }$T_{\Omega ,\alpha }$ and $M_{\Omega
,\alpha }$ {with a }$BMO$ functions {on vanishing generalized weighted
Morrey spaces, respectively. That is, in this paper we will consider this
problem.}

\section{Main results}

Let us state our main result as follows.

\begin{theorem}
\label{teo3}Suppose that $0<\alpha <n$, $1\leq s^{\prime }<p<\frac{n}{\alpha 
}$, $\frac{1}{q}=\frac{1}{p}-\frac{\alpha }{n}$, $1<q<\infty $, $\Omega \in
L_{s}(S^{n-1})\left( s>1\right) $ satisfies (\ref{f}) such that $k\in 
\mathbb{N}
$, $w\left( x\right) ^{s^{\prime }}\in A\left( \frac{p}{s^{\prime }},\frac{q%
}{s^{\prime }}\right) $, $A\in BMO\left( {\mathbb{R}^{n}}\right) $, $%
T_{\Omega ,\alpha }^{A,k}$, $M_{\Omega ,\alpha }^{A,k}$ are defined as (\ref%
{f1}), (\ref{f2}) and $T_{\Omega ,\alpha }^{A,k}$ satisfies (13) in \cite%
{Gurbuz3}. If $\varphi _{1}\in \mathcal{B}\left( w^{p}\right) $, $\varphi
_{2}\in \mathcal{B}\left( w^{q}\right) $ and the pair $(\varphi _{1},\varphi
_{2})$ satisfies the conditions 
\begin{equation}
c_{\delta }:=\dint\limits_{\delta }^{\infty }\left( 1+\ln \frac{t}{r}\right)
^{k}\sup_{x\in {\mathbb{R}^{n}}}\frac{\varphi _{1}\left( x,t\right) }{\left(
w^{q}\left( B\left( x,t\right) \right) \right) ^{\frac{1}{q}}}\frac{1}{t}%
dt<\infty  \label{6}
\end{equation}%
for every $\delta >0$, and 
\begin{equation}
\int\limits_{r}^{\infty }\left( 1+\ln \frac{t}{r}\right) ^{k}\frac{\varphi
_{1}\left( x,t\right) }{\left( w^{q}\left( B\left( x,t\right) \right)
\right) ^{\frac{1}{q}}}\frac{1}{t}dt\lesssim \frac{\varphi _{2}(x,r)}{\left(
w^{q}\left( B\left( x,t\right) \right) \right) ^{\frac{1}{q}}},  \label{7}
\end{equation}%
then the operator $T_{\Omega ,\alpha }^{A,k}$ is bounded from $VM_{p,\varphi
_{1}}\left( w^{p}\right) $ to $VM_{q,\varphi _{2}}\left( w^{q}\right) $.
Moreover,%
\begin{equation}
\left\Vert T_{\Omega ,\alpha }^{A,k}f\right\Vert _{VM_{q,\varphi _{2}}\left(
w^{q},{\mathbb{R}^{n}}\right) }\lesssim \left\Vert A\right\Vert _{\ast
}^{k}\left\Vert f\right\Vert _{VM_{p,\varphi _{1}}\left( w^{p},{\mathbb{R}%
^{n}}\right) },  \label{8}
\end{equation}%
\begin{equation*}
\left\Vert M_{\Omega ,\alpha }^{A,k}f\right\Vert _{VM_{q,\varphi _{2}}\left(
w^{q},{\mathbb{R}^{n}}\right) }\lesssim \left\Vert A\right\Vert _{\ast
}^{k}\left\Vert f\right\Vert _{VM_{p,\varphi _{1}}\left( w^{p},{\mathbb{R}%
^{n}}\right) }.
\end{equation*}
\end{theorem}

For $\alpha =0$, from Theorem \ref{teo3}, we get the following:

\begin{corollary}
Suppose that $1<p<\infty $, $s^{\prime }<p$, $\Omega \in
L_{s}(S^{n-1})\left( s>1\right) $ satisfies (\ref{f}) and (\ref{0*}) such
that $k\in 
\mathbb{N}
$, $w\left( x\right) ^{s^{\prime }}\in A_{\frac{p}{s^{\prime }}}$, $A\in
BMO\left( {\mathbb{R}^{n}}\right) $, $T_{\Omega }^{A,k}$, $M_{\Omega }^{A,k}$
are defined as%
\begin{eqnarray*}
T_{\Omega }^{A,k}f(x) &=&T_{\Omega }\left( \left( A\left( x\right) -A\left(
\cdot \right) \right) ^{k}f\left( \cdot \right) \right) \left( x\right)
,\qquad k=0,1,2,\ldots , \\
&=&p.v.\int\limits_{{\mathbb{R}^{n}}}\frac{\Omega (x-y)}{|x-y|^{n}}\left(
A\left( x\right) -A\left( y\right) \right) ^{k}f(y)dy
\end{eqnarray*}%
and the corresponding higher order ($=k$-th order) commutator operator of $%
M_{\Omega }$:
\end{corollary}

\begin{eqnarray*}
M_{\Omega }^{A,k}f(x) &=&M_{\Omega }\left( \left( A\left( x\right) -A\left(
\cdot \right) \right) ^{k}f\left( \cdot \right) \right) \left( x\right)
,\qquad k=0,1,2,\ldots , \\
&=&\sup_{r>0}\frac{1}{r^{n}}\int\limits_{|x-y|<r}\left\vert \Omega \left(
x-y\right) \right\vert \left\vert A\left( x\right) -A\left( y\right)
\right\vert ^{k}\left\vert f(y)\right\vert dy
\end{eqnarray*}%
and $T_{\Omega }^{A,k}$ satisfies (11) in \cite{Gurbuz4}. If $\varphi \in 
\mathcal{B}\left( w\right) $ and the pair $(\varphi _{1},\varphi _{2})$
satisfies the conditions%
\begin{equation*}
c_{\delta ^{\prime }}:=\dint\limits_{\delta }^{\infty }\left( 1+\ln \frac{t}{%
r}\right) ^{k}\sup_{x\in {\mathbb{R}^{n}}}\frac{\varphi _{1}\left(
x,t\right) }{\left( w^{p}\left( B\left( x,t\right) \right) \right) ^{\frac{1%
}{p}}}\frac{1}{t}dt<\infty
\end{equation*}%
for every $\delta ^{\prime }>0$, and 
\begin{equation*}
\int\limits_{r}^{\infty }\left( 1+\ln \frac{t}{r}\right) ^{k}\frac{\varphi
_{1}\left( x,t\right) }{\left( w^{p}\left( B\left( x,t\right) \right)
\right) ^{\frac{1}{p}}}\frac{1}{t}dt\lesssim \frac{\varphi _{2}(x,r)}{\left(
w^{p}\left( B\left( x,t\right) \right) \right) ^{\frac{1}{p}}},
\end{equation*}%
Then,%
\begin{equation*}
\left\Vert T_{\Omega }^{A,k}f\right\Vert _{VM_{p,\varphi _{2}}\left( w,{%
\mathbb{R}^{n}}\right) }\lesssim \left\Vert A\right\Vert _{\ast
}^{k}\left\Vert f\right\Vert _{VM_{p,\varphi _{1}}\left( w,{\mathbb{R}^{n}}%
\right) },
\end{equation*}%
\begin{equation*}
\left\Vert M_{\Omega }^{A,k}f\right\Vert _{VM_{p,\varphi _{2}}\left( w,{%
\mathbb{R}^{n}}\right) }\lesssim \left\Vert A\right\Vert _{\ast
}^{k}\left\Vert f\right\Vert _{VM_{p,\varphi _{1}}\left( w,{\mathbb{R}^{n}}%
\right) }.
\end{equation*}

\section{Proof of the main result}

\textbf{Proof of Theorem \ref{teo3}. }

\begin{proof}
By Definition \ref{definition1}, (13) in \cite{Gurbuz3} and (\ref{7}) we get%
\begin{eqnarray*}
\left\Vert T_{\Omega ,\alpha }^{A,k}f\right\Vert _{VM_{q,\varphi _{2}}\left(
w^{q},{\mathbb{R}^{n}}\right) } &=&\sup\limits_{x\in {\mathbb{R}^{n},r>0}}%
\frac{\left\Vert T_{\Omega ,\alpha }^{A,k}f\right\Vert _{L_{q}\left(
w^{q},B(x,r)\right) }}{\varphi _{2}(x,r)} \\
&\lesssim &\sup\limits_{x\in {\mathbb{R}^{n},r>0}}\frac{1}{\varphi _{2}(x,r)}%
\left\Vert A\right\Vert _{\ast }^{k}\left( w^{q}\left( B(x,r)\right) \right)
^{\frac{1}{q}} \\
&&\times \dint\limits_{r}^{\infty }\left( 1+\ln \frac{t}{r}\right)
^{k}\left\Vert f\right\Vert _{L_{p}\left( w^{p},B(x,t)\right) }\left(
w^{q}\left( B(x,t)\right) \right) ^{-\frac{1}{q}}\frac{1}{t}dt \\
&\lesssim &\sup\limits_{x\in {\mathbb{R}^{n},r>0}}\frac{1}{\varphi _{2}(x,r)}%
\left\Vert A\right\Vert _{\ast }^{k}\left( w^{q}\left( B(x,r)\right) \right)
^{\frac{1}{q}} \\
&&\times \dint\limits_{r}^{\infty }\left( 1+\ln \frac{t}{r}\right) ^{k}\frac{%
\varphi _{1}\left( x,t\right) }{\left( w^{q}\left( B\left( x,t\right)
\right) \right) ^{\frac{1}{q}}}\frac{1}{t}\left\Vert f\right\Vert
_{L_{p}\left( w^{p},B(x,t)\right) }\varphi _{1}(x,t)^{-1}dt \\
&\lesssim &\left\Vert A\right\Vert _{\ast }^{k}\left\Vert f\right\Vert
_{VM_{p,\varphi _{1}}\left( w^{p},{\mathbb{R}^{n}}\right) }.
\end{eqnarray*}%
At last, we need to prove that%
\begin{equation*}
\lim\limits_{r\rightarrow 0}\sup\limits_{x\in {\mathbb{R}^{n}}}\frac{1}{%
\varphi _{2}(x,r)}\left\Vert T_{\Omega ,\alpha }^{A,k}f\right\Vert
_{L_{q}\left( w^{q},B(x_{0},r)\right) }\lesssim \lim\limits_{r\rightarrow
0}\sup\limits_{x\in {\mathbb{R}^{n}}}\frac{1}{\varphi _{1}(x,r)}\Vert f\Vert
_{L_{p}\left( w^{p},B\left( x_{0},r\right) \right) }=0.
\end{equation*}%
Indeed, for any $\epsilon >0$, let $0<r<\psi $. By (13) in \cite{Gurbuz3},
we have%
\begin{equation*}
\frac{\left\Vert T_{\Omega ,\alpha }^{A,k}f\right\Vert _{L_{q}\left(
w^{q},B(x,r)\right) }}{\varphi _{2}(x,r)}\lesssim \left[ \mathcal{F}_{\psi
}\left( x,r\right) +\mathcal{G}_{\psi }\left( x,r\right) \right] ,
\end{equation*}%
where%
\begin{equation*}
\mathcal{F}_{\psi }\left( x,r\right) :=\frac{\left\Vert A\right\Vert _{\ast
}^{k}\left( w^{q}\left( B(x,r)\right) \right) ^{\frac{1}{q}}}{\varphi
_{2}(x,r)}\dint\limits_{r}^{\psi }\left( 1+\ln \frac{t}{r}\right)
^{k}\left\Vert f\right\Vert _{L_{p}\left( w^{p},B(x,t)\right) }\left(
w^{q}\left( B(x,t)\right) \right) ^{-\frac{1}{q}}\frac{1}{t}dt
\end{equation*}%
and%
\begin{equation*}
\mathcal{G}_{\psi }\left( x,r\right) :=\frac{\left\Vert A\right\Vert _{\ast
}^{k}\left( w^{q}\left( B(x,r)\right) \right) ^{\frac{1}{q}}}{\varphi
_{2}(x,r)}\dint\limits_{\psi }^{\infty }\left( 1+\ln \frac{t}{r}\right)
^{k}\left\Vert f\right\Vert _{L_{p}\left( w^{p},B(x,t)\right) }\left(
w^{q}\left( B(x,t)\right) \right) ^{-\frac{1}{q}}\frac{1}{t}dt.
\end{equation*}%
For $\sup\limits_{x\in {\mathbb{R}^{n}}}\sup\limits_{0<r<t}\frac{\left\Vert
f\right\Vert _{L_{q}\left( w^{q},B(x,r)\right) }}{\varphi _{2}(x,r)}<\frac{%
\epsilon }{2}$, we can select any constant $\psi >0$. This allows to guess
the first term properly from the type $r\in \left( 0,\psi \right) $ such that%
\begin{equation*}
\sup\limits_{x\in {\mathbb{R}^{n}}}\mathcal{F}_{\psi }\left( x,r\right) <%
\frac{\epsilon }{2}.
\end{equation*}

For the second term, in view of (\ref{6}), we obtain%
\begin{equation*}
\mathcal{G}_{\psi }\left( x,r\right) \lesssim \left\Vert A\right\Vert _{\ast
}^{k}\left\Vert f\right\Vert _{VM_{p,\varphi }\left( w^{p},{\mathbb{R}^{n}}%
\right) }\frac{\left( w^{q}\left( B(x,r)\right) \right) ^{\frac{1}{q}}}{%
\varphi _{2}(x,r)}.
\end{equation*}%
Since $\varphi _{2}\in \mathcal{B}\left( w^{q}\right) $, it gets along to
select $r$ minor sufficient such that 
\begin{equation*}
\sup\limits_{x\in {\mathbb{R}^{n}}}\frac{w^{q}\left( B(x,r)\right) }{\varphi
_{2}(x,r)}\lesssim \left( \frac{\epsilon }{2\left\Vert A\right\Vert _{\ast
}^{k}\left\Vert f\right\Vert _{VM_{p,\varphi }\left( w^{p},{\mathbb{R}^{n}}%
\right) }}\right) ^{q}.
\end{equation*}%
Hence,%
\begin{equation*}
\sup\limits_{x\in {\mathbb{R}^{n}}}\mathcal{F}_{\psi }\left( x,r\right) <%
\frac{\epsilon }{2}.
\end{equation*}%
Thus,%
\begin{equation*}
\frac{\left\Vert T_{\Omega ,\alpha }^{A,k}f\right\Vert _{L_{q}\left(
w^{q},B(x,r)\right) }}{\varphi _{2}(x,r)}<\epsilon .
\end{equation*}%
Therefore,%
\begin{equation*}
\lim\limits_{r\rightarrow 0}\sup\limits_{x\in {\mathbb{R}^{n}}}\frac{1}{%
\varphi _{2}(x,r)}\left\Vert T_{\Omega ,\alpha }^{A,k}f\right\Vert
_{L_{q}\left( w^{q},B(x_{0},r)\right) }=0.
\end{equation*}%
As a result, (\ref{8}) holds. On the other hand, since $M_{\Omega ,\alpha
}^{A,k}f(x)\leq \widetilde{T}_{\left\vert \Omega \right\vert ,\alpha
}^{A,k}\left( \left\vert f\right\vert \right) (x)$, $x\in {\mathbb{R}^{n}}$
(see Lemma 6 in \cite{Gurbuz3}) we can also use the same method for $%
M_{\Omega ,\alpha }^{A,k}$, so we omit the details. As a result, we complete
the proof of Theorem \ref{teo3}.
\end{proof}

\end{document}